\documentclass[onecolumn,draftcls]{IEEEtran}

\usepackage{amsmath,graphicx,bm,amssymb}

\renewcommand{\(}{\left(}
\renewcommand{\)}{\right)}
\renewcommand{\[}{\left[}
\renewcommand{\]}{\right]}
\renewcommand{\c}{\mathbf{c}}

\renewcommand{\S}{\mathbf{S}}

\newcommand{\thet}{\bm{\theta}}
\newcommand{\Nabla}{\bm{\nabla}}

\newcommand{\Del}{\bm{\Delta}}

\renewcommand{\H }{\mathbf{H}}

\newcommand{\0}{\mathbf{0}}
\newcommand{\D}{\mathbf{D}}

\newcommand{\x}{\mathbf{x}}
\newcommand{\I}{\mathbf{I}}

\newcommand{\A}{\mathbf{A}}

\newcommand{\U}{\mathbf{U}}

\newcommand{\Q}{\mathbf{Q}}
\newcommand{\K}{\mathbf{K}}

\newcommand{\X}{\mathbf{X}}

\newcommand{\B}{\mathbf{B}}

\newcommand{\Tr}[1]{{\rm{Tr}}\left\{#1\right\}}
\newcommand{\Trr}[1]{{\rm{Tr}}^2\left\{#1\right\}}
\renewcommand{\log}[1]{{\rm{log}}#1}
\renewcommand{\arg}[1]{{\rm{arg}}#1}

\newtheorem{theorem}{Theorem}

\begin{document}

\title{Covariance estimation in decomposable Gaussian graphical models}

\author{Ami Wiesel, Yonina C. Eldar and Alfred O. Hero III \thanks{A. Wiesel and A. O. Hero are with the Department of Electrical Engineering and Computer Science, University of Michigan, Ann Arbor, MI 48109, USA. Tel: (734)763-0564, Fax: (734) 763-8041. Emails: \{amiw,hero\}@umich.edu.} \thanks{Y.C. Eldar is with the Technion - Israel Institute of Technology, Haifa, Israel 32000. Email: yonina@ee.technion.ac.il.} \thanks{The work of A. Wiesel was supported by a Marie Curie Outgoing International Fellowship within the 7th European Community Framework Programme.} }

\maketitle

\begin{abstract}
  Graphical models are a framework for representing and exploiting prior conditional independence structures within distributions using graphs. In the Gaussian case, these models are directly related to the sparsity of the inverse covariance (concentration) matrix and allow for improved covariance estimation with lower computational complexity. We consider concentration estimation with the mean-squared error (MSE) as the objective, in a special type of model known as decomposable. This model includes, for example, the well known banded structure and other cases encountered in practice. Our
  first contribution is the derivation and analysis of the minimum variance unbiased estimator (MVUE) in decomposable graphical models. We provide a simple closed form solution to the MVUE and compare it with the classical maximum likelihood estimator (MLE) in terms of performance and complexity. Next, we extend the celebrated Stein's unbiased risk estimate (SURE) to graphical models. Using SURE, we prove that the MSE of the MVUE is always smaller or equal to that of the biased MLE, and that the MVUE itself is dominated by other approaches.  In addition, we propose the use of SURE as a constructive mechanism for deriving new covariance estimators. Similarly to the classical MLE, all of our proposed estimators have simple closed form solutions but result in a significant reduction in MSE.
\end{abstract}

\section{Introduction}

Covariance estimation in Gaussian distributions is a classical and
fundamental problem in statistical signal processing.  Many
applications, varying from array processing to functional
genomics, rely on accurately estimated covariance matrices
\cite{krim:96,dougherty:05}. Recent interest in inference in high
dimensional settings using small sample sizes has caused the topic
to rise to prominence once again. A natural approach in these
settings is to incorporate additional prior knowledge in the form
of structure and/or sparsity in order to ensure stable estimation.
Gaussian graphical models provide a method of representing
conditional independence structure among the different variables
using graphs. An important property of the Gaussian distributions
is that conditional independence among groups of variables is
associated with sparsity in the inverse covariance. Due to the sparsity, these models allow for
efficient implementation of statistical inference algorithms,
e.g., the iterative proportional scaling technique \cite{lauritzen:book,dempster:72}.

Over the last years, statistical graphical models were successfully applied to speech recognition \cite{jbilmes:2000,jbilmes:2005}, image processing \cite{willsky02multiresolution,choiWillsky:07}, sensor networks \cite{willsky:SPM}, computer networks \cite{wiesel:dpca} and other fields in signal processing. Efficient Bayesian inference in Gaussian graphical models is well established \cite{sudderth:2004,mali:2006,chandras:2008}. Estimation of deterministic parameters received less attention, but have also been considered implicitly, e.g., the recent works on inverse covariance structure in the context of state of the art array processing \cite{Abramovich:2007,Abramovich:nov2008ii}.

Estimation of deterministic parameters in Gaussian graphical
models is basically covariance estimation since the Gaussian
distribution is completely parameterized by second order
statistics. The most common approach to covariance estimation is
maximum likelihood. When no prior information is available, this
method yields the sample covariance matrix. It is asymptotically
unbiased and efficient but does not minimize the mean-squared
error (MSE) in general. Indeed, depending on the performance
measure, better estimators can be obtained through regularization,
shrinkage, empirical Bayes and other methods
\cite{stein:rietz,haff:empirical,Tsukuma:2006,ledoitwolf:04,yang:94,bickel:2008,bickel:2008th}.

Covariance estimation in Gaussian graphical models involves the estimation of the unknown covariance based on the observed realizations and prior knowledge of the conditional independence structure within the distribution \cite{dempster:72,lauritzen:book,speed:86,dawid:93}. The prior information allows for better performance with lower computational complexity. Decomposable graphical models, also known as chordal or triangulated, satisfy a special structure which leads to a simple closed form solution to the maximum likelihood estimate (MLE) as well as elegant analysis. These models include many practical signal processing structures such as the banded concentration matrix and its variants \cite{Abramovich:2007,Kavcic:2000,bickel:2008,bickel:2008th} as well as multiscale structures \cite{willsky02multiresolution,choiWillsky:07}.

Covariance selection is a related topic which addresses the joint problem of covariance estimation and graphical model selection. This setting is suitable to many modern applications in which the conditional independence structure is unknown and must be learned from the observations. Numerous selection methods have been recently considered for both arbitrary graphical models \cite{friedman-2007,yuan:2007,banerjee-2007,anandkumar:2009} and decomposable models\cite{Deshpande:2001,Jones05experimentsin}. Clearly, these methods are intertwined with covariance estimation. For example, the latter is the core of the intermediate steps in many of the advanced greedy stepwise selection methods.

In this paper, we consider inverse covariance estimation in
decomposable Gaussian graphical models with the mean-squared error
(MSE) as our objective. Except for the prior conditional
independence structure, we do not assume any prior knowledge on
the covariance and treat it as an unknown deterministic parameter.
Our main contribution is the derivation of the minimum variance
unbiased estimator (MVUE) of the inverse covariance. Similarly to
the MLE, the MVUE has a simple closed form solution which can be
efficiently implemented in a distributed manner. Moreover, it
minimizes the MSE among all unbiased estimators. We also prove
that it has smaller MSE than the biasd MLE. The proof is based on
an extension of the celebrated Stein unbiased risk estimate (SURE)
\cite{stein:rietz,Hudson:1978,eldar:gsure} to Gaussian graphical
models. Using SURE we prove that the MVUE dominates the MLE in
terms of MSE, i.e., its MSE is always smaller or equal to that of the MLE. In addition, we prove that the MVUE itself is dominated by other
biased estimators. Next, we propose the use of SURE as a method for
hyper-parameter tuning in existing covariance estimation approaches, e.g., the
conjugate prior based methods proposed in \cite{letac:2007,rajaratman:07}.

The outline of the paper is as follows. We begin in Section \ref{sec_models} where we formally define decomposable graphical models, provide a few illustrative applications, and formulate the estimation problem. In Section \ref{sec_MLE}, we review the classical MLE approach and derive the finite sample MVUE. Next, in Section \ref{sec_sure} we consider SURE and its applications. While our estimators have lower MSE, they require more samples in order to ensure positive semidefiniteness. This issue is addressed in Section \ref{sec_pospart}. We evaluate the performance of the different estimators using numerical simulations in Section \ref{sec_num}, and offer concluding remarks in Section \ref{sec_conc}.

The following notation is used. Boldface upper case letters denote
matrices, boldface lower case letters denote column vectors, and
standard lower case letters denote scalars. The superscripts
$(\cdot)^T$ and $(\cdot)^{-1}$ denote the transpose and matrix
inverse, respectively. For sets $a$ and $b$, the cardinality is
denoted by $|a|$ and the set difference operator is denoted by
$a\setminus b$.  The operator $\|\X\|$ denotes the Frobenius norm
of a matrix $\X$, namely $\|\X\|^2=\operatorname{Tr}(\X^T\X)$, and
$\X\succ\0$ means that $\X$ is positive definite. The zero fill-in
operator $\[\cdot\]^0$ outputs a conformable matrix where the
argument occupies its appropriate sub-block and the rest of the
matrix has zero valued elements (See \cite{lauritzen:book} for the
exact definition of this operator).

\section{Covariance estimation in graphical models}\label{sec_models}
In this section, we provide an introduction to decomposable Gaussian graphical models based on \cite{lauritzen:book} along with a few motivating applications for their use in modern statistical signal processing. Then, we formulate the inverse covariance estimation problem addressed in this paper.

\subsection{Decomposable Gaussian graphical models}
Graphical models are intuitive characterizations of conditional independence structures within distributions. An undirected graph $\mathcal{G}=\(V,E\)$ is a set of nodes $V=\{1,\cdots,|V|\}$ connected by undirected edges $E=\{\(i_1,j_1\),\cdots\(i_{|E|},j_{|E|}\)\}$, where we use the convention that each node is connected to itself, i.e., $\(i,i\)\in E$ for all $i\in V$. Let $\x$ be a zero mean random vector of length $p=|V|$ whose elements are indexed by the nodes in $V$. The vector $\x$ satisfies the Markov property with respect to $\mathcal{G}$, if for any pair of non-adjacent nodes the corresponding pair of elements in $\x$ are conditionally independent of the remaining elements, i.e.,  $\[\x\]_i$ and $\[\x\]_j$ are conditionally independent of $\[\x\]_{V\setminus i,j}$ for any $\{i,j\}\notin E$. For the multivariate Gaussian distribution, conditional indpendence is equivalent to conditional uncorrelation:
\begin{eqnarray}\label{uncorr}
  E\left\{\left.\(\[\x\]_i-E\left\{\left.\[\x\]_i\right| \[\x\]_{V\setminus i,j}\right\}\)
  \(\[\x\]_j-E\left\{\[\x\]_j\left| \[\x\]_{V\setminus i,j}\right.\right\}\)\right| \[\x\]_{V\setminus i,j}\right\}=0,\qquad \{i,j\}\notin E.
\end{eqnarray}
Simple algebraic manipulations (see Appendix \ref{app_sparsity}) show that (\ref{uncorr}) is equivalent to
\begin{eqnarray}
  \[\K\]_{i,j}=0\qquad {\text{for all}} \qquad \{i,j\}\notin E.
\end{eqnarray}
where $\K$ is the concentration (inverse covariance) matrix of $\x$ defined as
\begin{eqnarray}
  \K=\(E\left\{\x\x^T\right\}\)^{-1}.
\end{eqnarray}
To summarize,
Gaussian graphical models are directly related to sparsity in the
concentration matrix $\K$. Therefore, throughout this paper we
focus on $\K$ and parameterize the Gaussian distribution using the
concentration matrix rather than the covariance matrix.

Decomposable models are a special case of graphical models in
which the conditional independence graphs satisfy an appealing
structure. A graph (or subgraph) is complete if all of its nodes
are connected by an edge. A subset $c\subseteq V$ separates
$a\subseteq V$ and $b\subseteq V$ if all paths from $a$ to $b$
intersect $c$. A triple $\{a,c,b\}$ of disjoint subsets of
$\mathcal{G}$ is a decomposition of $\mathcal{G}$ if $V= a\cup
c\cup b$, $c$ separates $a$ and $b$, and $c$ is complete. Finally,
a graph $\mathcal{G}$ is decomposable if it is complete, or if it
can be recursively decomposed into two decomposable subgraphs
$\mathcal{G}_{a\cup c}$ and $\mathcal{G}_{c\cup b}$.

It is convenient to
represent a decomposable graph using cliques. A clique is a
maximal complete subset of nodes. If $\mathcal{G}$ is decomposable
then there exists a sequence of cliques $C_1,\cdots,C_K$ which
satisfy a {\em{perfect elimination order}}. This order implies
that $\{H_{j-1}\backslash S_j,S_j,R_j\}$ is a decomposition of the
subgraph $\mathcal{G}_{H_j}$ where
\begin{eqnarray}\label{HSR}
  H_{j}&=&C_1\cup C_2\cup\cdots\cup C_{j},\quad j=1,\cdots,K\\
  S_{j}&=&H_{j-1}\cap C_j, \quad j=2,\cdots,K\\
  R_{j}&=&H_j\backslash H_{j-1}, \quad j=2,\cdots,K.
\end{eqnarray}
For later use, we denote the cardinalities of the
cliques and separators by $|C_k|=c_k$ and $|S_k|=s_k$, so that
\begin{eqnarray}
  \sum_{k=1}^Kc_k-\sum_{k=2}^Ks_k=p.
\end{eqnarray}
Consequently, if a Gaussian multivariate $\x$ satisfies a decomposable graphical model over the graph $\mathcal{G}$, then its concentration matrix $\K$ belongs to the set of decomposable sparse positive semidefinite matrices:
\begin{eqnarray}\label{KK}
  {\mathcal{K}}&=&\left\{\K: \K\succeq\0, \qquad \[\K\]_{H_{j-1}\setminus S_j,R_j}=\0,\quad j=2,\cdots,K\right\}.
\end{eqnarray}

Decomposable graphical models appear in many signal processing
applications. We now review a few representative examples:
\begin{itemize}
  \item {\em{Diagonal or block diagonal:}} A trivial graphical model is the diagonal or block diagonal model, in which the cliques are non-overlapping. For example, the following matrix has two cliques $C_1=\{1,2\}$ and $C_2=\{3,4,5\}$:
      \begin{eqnarray}
        \includegraphics[width=0.15\textwidth]{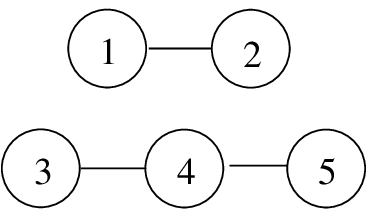}\qquad\qquad\qquad
        \K=\[\begin{array}{ccccc}
                 \Box & \Box &  &  &  \\
                 \Box & \Box &  &  &  \\
                  &  & \Box & \Box & \Box \\
                  &  & \Box & \Box & \Box \\
                  &  & \Box & \Box & \Box
               \end{array}\].
      \end{eqnarray}

  \item {\em{Two coupled blocks:}} The simplest non-trivial decomposable graphical model is the two coupled blocks. For example, the following matrix has two cliques $C_1=\{1,2,3\}$ and $C_2=\{3,4,5\}$ coupled through $S_2=\{3\}$:
      \begin{eqnarray}
      \includegraphics[width=0.15\textwidth]{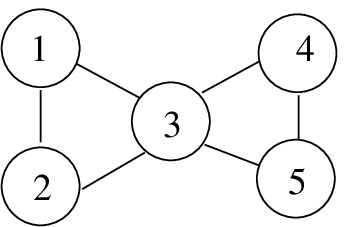}\qquad\qquad\qquad
        \K=\[\begin{array}{ccccc}
                 \Box & \Box & \Box &  &  \\
                 \Box & \Box & \Box &  &  \\
                 \Box & \Box & \Box & \Box & \Box \\
                  &  & \Box & \Box & \Box \\
                  &  & \Box & \Box & \Box
               \end{array}\].
      \end{eqnarray}
      Note that the definition of decomposable model is that they can be recursively subdivided into this main building block.
  \item {\em{Banding:}} Another frequently used decomposable graphical model is the $L$'th order banded structure in which only the $L+1$ principal diagonals of $\K$ have non-zero elements. For example, the following matrix is banded with $L=2$:
            \begin{eqnarray}
        \includegraphics[width=0.30\textwidth]{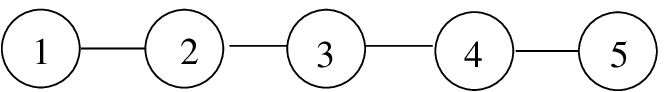}\qquad\qquad\qquad
        \K=\[\begin{array}{ccccc}
                 \Box & \Box &  &  &  \\
                 \Box & \Box & \Box  &  &  \\
                  &  \Box & \Box & \Box & \\
                  &  & \Box & \Box & \Box \\
                  &  &  & \Box & \Box
               \end{array}\],
      \end{eqnarray}
      and has for cliques $C_1=\{1,2\}$, $C_2=\{2,3\}$, $C_3=\{3,4\}$ and $C_4=\{4,5\}$.  It is appropriate whenever the indices of the multivariate represent physical quantities such as time or space, and the underlying assumption is that distant variables are conditionally independent of closer variables. A special case of this structure is the autoregressive (AR) model. The AR model is stationary and leads to a banded Toeplitz matrix. The more general banded graphical model corresponds to a non-stationary autoregressive process. It was recently shown that this structure is a good model for state-of-the-art radar systems \cite{Abramovich:2007} (see also \cite{Kavcic:2000}). A natural extension of the $L$'th banded model is differential banding in which multiple band lengths are utilized. It is straightforward to show that the corresponding graph is still decomposable with cliques of different cardinalities. This form was empirically validated to be a reasonable model in call center management in operations research \cite{rajaratman:07}.
  \item {\em{Arrow (Star):}} Another common decomposable graphical models takes the form of an arrow motif in the concentration matrix. This structure is appropriate when there is a single common global sub-block and numerous local sub-blocks which are conditionally independent given the global variables. For example, the following concentration matrix specifies an arrow graphical model
      \begin{eqnarray}
        \includegraphics[width=0.15\textwidth]{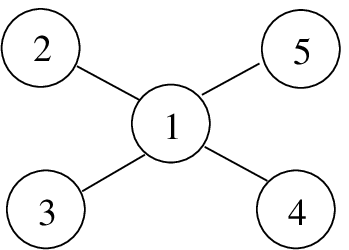}\qquad\qquad\qquad
        \K=\[\begin{array}{ccccc}
                 \Box & \Box & \Box & \Box & \Box \\
                 \Box & \Box &   &  &  \\
                 \Box &   & \Box & & \\
                 \Box &  &  & \Box & \\
                 \Box &  &  &  & \Box
               \end{array}\],
      \end{eqnarray}
      with cliques $C_1=\{1,2\}$, $C_2=\{1,3\}$, $C_3=\{1,4\}$ and $C_4=\{1,5\}$.  A typical signal processing application is a wireless network in which the global node is the access point and the local nodes are the terminals. Other applications of these models are discussed in \cite{sunsun:2005}.
  \item {\em{Multiscale:}} A common graphical model in image processing is based on the multiscale (multiresolution) framework. Here, the decomposable graph is a tree of nodes (or cliques) \cite{willsky02multiresolution,choiWillsky:07}:
      \begin{eqnarray}
        \includegraphics[width=0.25\textwidth]{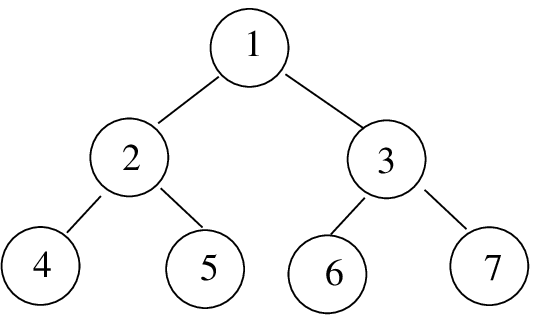}\qquad\qquad\qquad
        \K=\[\begin{array}{ccccccc}
                 \Box & \Box & \Box & & & &\\
                 \Box & \Box &  & \Box & \Box  &  & \\
                 \Box & & \Box & & & \Box & \Box \\
                  &  \Box & & \Box & & & \\
                  & \Box  & &  & \Box & & \\
                  & & \Box  &  &  & \Box & \\
                  & & \Box  &  &  & & \Box\\
               \end{array}\].
      \end{eqnarray}

\end{itemize}

\subsection{Problem formulation}\label{sec_prob}
We are now ready to state the problem addressed in this paper. Let
$\x$ be a length-$p$ zero mean Gaussian random vector, with
concentration matrix $\K\in{\mathcal{K}}$ as in (\ref{KK}). Given
$n$ independent realizations of $\x$ denoted by
$\{\x[i]\}_{i=1}^n$, and knowledge of the conditional independence
structure, our goal is to derive an estimate $\hat\K$ of $\K$ with
minimum MSE, where the MSE is defined as
\begin{eqnarray}\label{mse}
  {\text{MSE}}\(\K\)&=&E\left\{\left\|\hat\K-\K\right\|^2\right\}.
\end{eqnarray}
Here the norm is the matrix Frobenius norm.  The MSE in
(\ref{mse}) is a function of the unknown parameter $\K$ and cannot
be minimized directly. This dependency is the main difficulty in
minimum MSE estimation of deterministic parameters, in contrast to
the Bayesian framework in which the MSE is a function of the
distribution of $\K$ but not of $\K$ itself. More details on this issue can be found in \cite{kay2008rbe}.

Due to the difficulty of minimum MSE estimation, it is customary to restrict attention to unbiased estimators. For this purpose, the MSE is decomposed into its squared bias and variance components defined as
\begin{eqnarray}
  {\text{MSE}}\(\K\)&=&{\text{BIAS}}^2\(\K\)+{\text{VAR}}\(\K\)
\end{eqnarray}
where
\begin{eqnarray}
  {\text{BIAS}}^2\(\K\)&=&\left\|E\left\{\hat\K\right\}-\K\right\|^2\nonumber\\
  {\text{VAR}}\(\K\)&=&E\left\{\left\|\hat\K-E\left\{\hat\K\right\}\right\|^2\right\}.
\end{eqnarray}
We call $\hat\K$ an unbiased estimator if ${\text{BIAS}}\(\K\)=0$.
Although the variance may also depend on $\K$, asymptotically in
many cases an unbiased estimate exists that minimizes the MSE.

Our choice of MSE of $\K$ as a performance measure requires
further elaboration. There are numerous competing metrics which
could have been adopted:  MSE of $\K^{-1}$; weighted norms;
Kullbuck-Leibler based distances and others. Each of these
measures will lead to different estimators. We focus on the MSE of
the inverse covariance due to the following reasons. Graphical
models specify the structure of the concentration matrix rather
than the structure of the covariance matrix so that the
ocnetnraiton seems to be a more intuitive paramter. Furthermore,
the concentration is the natural parameter of the Gaussian
distribution as an exponential family \cite{lauritzen:book}. It is
parameterized by the free variables associated with the cliques,
and has zero values elsewhere. In contrast, the covariance matrix
is not a natural parameter and has unnecessary and functionally
dependent variables outside the cliques. When included in the
performance measure, these variables mash the behavior of
the free variables within the cliques. Some of the results in this
paper, such as the SURE identity, can also be applied to other
performance measures.

\section{Maximum likelihood and minimum variance unbiased estimation}\label{sec_MLE}
In this section, we review the classical MLE approach to inverse covariance
estimation in decomposable Gaussian graphical models and then derive the MVUE estimator.

\subsection{Maximum likelihood estimation}
We begin with a short review of the sufficient statistics and MLE approach. For a more
detailed treatment the reader is referred to
\cite{lauritzen:book}.

When no prior information is available and the model consists of one clique $C_1=\{1,\cdots,p\}$, the model is said to be saturated. In this case, a minimal sufficient statistic for estimating $\K$ is the sample covariance matrix
\begin{eqnarray}
  \overline{\S}=\sum_{i=1}^n\x[i]\x^T[i].
\end{eqnarray}
Its distribution is Wishart with $n$ degrees of freedom and natural parameter $\K$
\begin{eqnarray}\label{wishart}
  p\(\overline{\S};\K\)&=&{\mathcal{W}}_p\(\overline{\S};n,\K\)\nonumber\\
&=&\frac{\left|\overline{\S}\right|^{\frac{n-p-1}{2}}\left|\K\right|^{\frac{n}{2}}}{2^{\frac{np}{2}}\Gamma_p\(\frac{n}{2}\)}e^{-\frac{1}{2}\Tr{\K\overline{\S}}}
\end{eqnarray}
where $\Gamma_p\(\cdot\)$ is the multivariate Gamma function. In
graphical models, the known conditional independence structure
specifies some of the entries of $\K$ and the complete sample
covariance is no longer necessary. A
minimal sufficient statistic is the incomplete matrix $\S$ defined
as
\begin{eqnarray}\label{defS}
  \[\S\]_{i,j}=\left\{\begin{array}{cc}
                                   \[\overline{\S}\]_{i,j} & \{i,j\}\in E \\
                                   ? & {\text{else}}
                                 \end{array}
  \right.
\end{eqnarray}
where $?$ denotes unspecified elements. In particular, if ${\mathcal{G}}$ is decomposable then only the sub-blocks associated with cliques and separators
\begin{eqnarray}
  \S_{{k}}&\equiv&\[\S\]_{C_k,C_k},\quad k=1,\cdots,K\nonumber\\
  \S_{[k]}&\equiv&\[\S\]_{S_k,S_k},\quad k=2,\cdots,K
  \end{eqnarray}
are necessary.
Their marginal distributions are also Wishart
\begin{eqnarray}
  p\(\S_{{k}};\(\[\K^{-1}\]_{C_k,C_k}\)^{-1}\)&=&{\mathcal{W}}_{c_k}\(\S_{k};n,\(\[\K^{-1}\]_{C_k,C_k}\)^{-1}\)\nonumber\\
  p\(\S_{[k]};\(\[\K^{-1}\]_{S_k,S_k}\)^{-1}\)&=&{\mathcal{W}}_{s_k}\(\S_{[k]};n,\(\[\K^{-1}\]_{S_k,S_k}\)^{-1}\).
\end{eqnarray}
Their joint distribution, denoted by $p\(\S\)$, is derived by marginalization of the distribution of $\overline{\S}$ over the unnecessary elements, and is known as the hyper Wishart distribution \cite{dawid:93}
\begin{eqnarray}\label{hyper}
  p\(\S;\K\)&=& \frac{\prod_kp\(\S_{k};\(\[\K^{-1}\]_{C_k,C_k}\)^{-1}\)}{\prod_kp\(\S_{[k]};\(\[\K^{-1}\]_{S_k,S_k}\)^{-1}\)}.
\end{eqnarray}
The simple product form of the distribution is a direct consequence of the decomposable structure.

Gaussian graphical models are regular models in the exponential
family. By appropriate definitions that transform the variables
from matrices to vectors and take into account the symmetry, both
the Wishart and the hyper Wishart distribution can be written as
natural exponential distributions. In Appendix \ref{app exp}, we
provide a few known results on this family that will be used in
the sequel.

The MLE of $\K$ is defined as\footnote{An alternative but equivalent definition is $\hat\K_{{\text{MLE}}}=\arg\max_{\K\in{\mathcal{K}}}\log p\(\overline{\S};\K\)$.}
\begin{eqnarray}
  \hat\K_{{\text{MLE}}}&=&\arg\max_{\K\succ\0}\log p\(\S;\K\),
\end{eqnarray}
and has the following closed form solution \cite{lauritzen:book}
\begin{eqnarray}\label{MLEst}
    \hat\K_{{\text{MLE}}}&=&\sum_{k=1}^K\[n\S_{{k}}^{-1}\]^0-\sum_{k=2}^K\[n\S_{[k]}^{-1}\]^0.
\end{eqnarray}
It exists with probability one if and only if $n\geq\max_k {c_k}$, in which case it is positive semidefinite. It is locally consistent in the sense that  the local and global versions of the cliques agree with each other:
\begin{eqnarray}\label{localconst}
   \[\hat\K_{\text{MLE}}^{-1}\]_{C_k,C_k}=\frac{1}{n}\S_{{k}},\qquad k=1,\cdots,K.
\end{eqnarray}
Both of these properties suggest that the MLE in a decomposable model performs as if the model was block diagonal with non-overlapping cliques $C_k$.

In general, the MLE is a biased estimator and does not minimize
the MSE. One of the main motivations for the MLE is that
asymptotically in $n$ it is an MVUE. Therefore, we now
address the finite sample MVUE in decomposable graphical models.
Interestingly, we will show that the MVUE does not behave as if
the model was block diagonal and improves performance by taking
into account the coupling between the cliques. We will also see
that it dominates the MLE, namely its MSE is smaller for all
possible values of $\K$.

\subsection{Minimum variance unbiased estimation}\label{sec_mvue}
For finite sample size, the MVUE is provided in the following theorem.
\begin{theorem}\label{theorem_mvue}
  The estimator
 \begin{eqnarray}\label{Kub}
    \hat\K_{\text{MVUE}}&=&\sum_{k=1}^K\[\(n-{c_k}-1\)\S_{{k}}^{-1}\]^0-\sum_{k=2}^K\[\(n-{s_k}-1\)\S_{[k]}^{-1}\]^0
\end{eqnarray}
exists with probability one if $n\geq\max_k {c_k}$, and is the minimum variance unbiased estimator (MVUE) of $\K\in{\mathcal{K}}$ given the incomplete sample covariance matrix $\S$.
\end{theorem}

The proof is available in Appendix \ref{app exp} and is based on the general MVUE for exponential family distributions.


Theorem 1 specifies the MVUE of $\K$ in decomposable graphical models. The estimator is similar in structure and complexity to the MLE in (\ref{MLEst}) and it is easy to see that asymptotically in $n$ the two estimators are equivalent. In the saturated case the MVUE is a scaled version of the MLE. In many signal processing applications (e.g., principal component analysis) the overall performance is indifferent to a change in scaling of the covariance. In decomposable graphical models, the MVUE is not a simple re-scaling of the MLE and there may be practical performance gains to its use with almost no additional cost in computational complexity.

Recall that the MLE requires only $n\geq\max_k {c_k}$ samples in
order for it to exist and be positive definiteness. This
is not true for the MVUE which may require more samples to ensure
positive semidefiniteness. For example, consider a simple $K=2$
cliques model. Using the matrix inversion formula for partitioned
matrices \cite[page 572]{kay:93}, it can be verified that
\begin{eqnarray}\label{invKubabc}
\[\hat\K_{\text{MVUE}}^{-1}\]_{S_2,S_2}=\frac{1}{n-p-1}\[\S\]_{S_2,S_2}.
\end{eqnarray}
A necessary condition for positive definiteness with probability
one of $\hat\K_{\text{MVUE}}$ is
$[\hat\K_{\text{MVUE}}^{-1}]_{S_2,S_2}\succ\0$ which is equivalent
to $n>p+1$. Thus, although $n\geq\max_k {c_k}$ suffices for existence
of an estimator it is meaningless unless $n>p+1$. Identity
(\ref{invKubabc}) may suggest that the MVUE is locally consistent
but it can be verified that this is not true, i.e.,
\begin{equation}
\[\hat\K_{\text{MVUE}}^{-1}\]_{C_k,C_k}\neq\frac{1}{n-p-1}\S_{{k}}
\end{equation}
for $k=1,2$. Evidently, in contrast to the MLE, the MVUE does not
behave as if the model were block diagonal and it accounts for the
coupling between the cliques.

The MVUE minimizes the MSE over the class of unbiased estimators. This is an important property but it does not ensure optimality over all estimators, biased or unbiased. In the next section, we prove that the MVUE actually dominates the biased MLE in terms of MSE performance.

\section{Stein's unbiased risk estimate (SURE)}\label{sec_sure}
Stein's unbiased risk estimate (SURE) provides an unbiased approximation of the MSE. The SURE approach was originally applied to the estimation of a Gaussian mean parameter \cite{stein:rietz}. It was generalized to the Wishart distribution in \cite{stein:rietz,haff:empirical}, and later extended to the natural parameters of any exponential family distribution in \cite{Hudson:1978,eldar:gsure}. The following theorem extends these results to Gaussian graphical models.

\begin{theorem}\label{theoremSURE}
  Let $\S$ be an incomplete sample covariance matrix associated with a graphical model $\mathcal{G}=\{V,E\}$ as defined in (\ref{defS}), and let $\H\(\S\)$ be a continuously differentiable matrix function of $\S$. Then,
  \begin{eqnarray}\label{identityeq}
  E\left\{\Tr{\H\(\S\)\K}\right\}=E\left\{\Tr{\H\(\S\)\hat\K_{\text{MVUE}}+2\Nabla\H\(\S\)}\right\}
  \end{eqnarray}
where $\hat\K_{\text{MVUE}}$ is the MVUE of $\K$ given $\S$ and
the differential operator is a $p\times p$ matrix with elements
\begin{eqnarray}\label{nabla}
  \Nabla=\left\{\begin{array}{ll}
                                \frac{\partial}{\partial \S_{ij}} & i=j \\
                                \frac{1}{2}\frac{\partial}{\partial \S_{ij}} & i\neq j,\{i,j\}\in E \\
                                0 & \text{else}.
                              \end{array}
  \right.
\end{eqnarray}
\end{theorem}

The proof is available in Appendix \ref{app exp} and is based on the general SURE approach to estimation in exponential family distributions. The theorem holds for any Gaussian graphical model but is useful only when the MVUE can be evaluated. In the decomposable case, the MVUE is provided in (\ref{Kub}) and the differential operator has the appealing form:
\begin{eqnarray}
  \Nabla&=&\sum_{k=1}^K\[\Nabla_k\]^0-\sum_{k=2}^K\[\Nabla_{[k]}\]^0
\end{eqnarray}
where $\Nabla_{k}=\[\Nabla\]_{C_k,C_k}$ and $\Nabla_{[k]}=\[\Nabla\]_{S_k,S_k}$ are the $c_k\times c_k$ and $s_k\times s_k$ differential operators within the saturated cliques and the saturated separators, respectively.

In the following subsections, we apply SURE to the derivation and MSE analysis of several estimators.

\subsection{MVUE dominates MLE}
Our first application of Theorem \ref{theoremSURE} is to
prove that the MLE is inadmissible and dominated by the MVUE:

\begin{theorem}
  The MVUE in (\ref{Kub}) dominates the MLE in (\ref{MLEst}) in terms of MSE:
\begin{eqnarray}\label{muve_ml}
  E\left\{\left\|\hat\K_{\text{MVUE}}-\K\right\|^2\right\}\leq E\left\{\left\|\hat\K_{\text{MLE}}-\K\right\|^2\right\},
\end{eqnarray}
for all $\K$ in the set ${\mathcal{K}}$ defined in (\ref{KK}).
\end{theorem}
\begin{proof}
  The difference in MSEs is
\begin{eqnarray}\label{deltaml}
  \delta&=&E\left\{\left\|\hat\K_{\text{MVUE}}-\K\right\|^2\right\}-E\left\{\left\|\hat\K_{\text{MLE}}-\K\right\|^2\right\}\nonumber\\
  &=&E\left\{\left\|\hat\K_{\text{MVUE}}\right\|^2-\left\|\hat\K_{\text{MLE}}\right\|^2-2\Tr{\H\K}\right\}\nonumber\\
  &=&E\left\{-\left\|\H \right\|^2-4\Tr{\Nabla\H}\right\}
\end{eqnarray}
where we applied Theorem \ref{theoremSURE} with
\begin{eqnarray}
  \H&=&\hat\K_{\text{MVUE}}-\hat\K_{\text{MLE}}\\
  &=&-\[\sum_{k=1}^K\({c_k}+1\)\[\S_{{k}}^{-1}\]^0-\sum_{k=2}^K\({s_k}+1\)\[\S_{[k]}^{-1}\]^0\].\nonumber
\end{eqnarray}
Therefore, in order to prove that $\delta\leq 0$  it is sufficient to show that
\begin{align}
\Tr{\Nabla\H }\geq 0.
\end{align}
From \cite[(5.4)iii]{haff:empirical}:
\begin{eqnarray}\label{NablaiS}
  \Tr{\Nabla_k\S_k^{-1}}&=&-\frac{1}{2}\Tr{\S_k^{-2}}-\frac{1}{2}\Trr{\S_k^{-1}}\nonumber\\
  \Tr{\Nabla_{[k]}\S_{[k]}^{-1}}&=&-\frac{1}{2}\Tr{\S_{[k]}^{-2}}-\frac{1}{2}\Trr{\S_{[k]}^{-1}}.
\end{eqnarray}
Therefore,
\begin{align}
&\Tr{\Nabla\H }=
\frac{1}{2}\[\sum_{k=1}^K\({c_k}+1\)\[\Tr{\S_{{k}}^{-2}}+\Trr{\S_{{k}}^{-1}}\]-\sum_{k=2}^K\({s_k}+1\)\[\Tr{\S_{[k]}^{-2}}
+\Trr{\S_{[k]}^{-1}}\]\]\nonumber\\
&\geq\frac{1}{2}\({s_k}+1\)\[\sum_{k=2}^K\Tr{\S_{{k}}^{-2}}-\Tr{\S_{[k]}^{-2}}+\Trr{\S_{{k}}^{-1}}
-\Trr{\S_{[k]}^{-1}}\]\geq 0
\end{align}
where we have used ${c_k}\geq {s_k}$ for $k=2,\cdots,K$ and
\begin{align}
\label{CvS}&\Tr{\S_{{k}}^{-1}}\geq\Tr{\S_{[k]}^{-1}},\quad k=2,\cdots,K,\\
\label{CvS2}&\Tr{\S_{{k}}^{-2}}\geq\Tr{\S_{[k]}^{-2}},\quad
k=2,\cdots,K.
\end{align}
The last two inequalites follow from Appendix \ref{app_proofs}.
\end{proof}

\subsection{MVUE is inadmissible}
We now use SURE to prove that the MVUE is inadmissible since its
MSE can be improved upon by another biased estimator.

\begin{theorem}\label{theorem_scaled}
  The estimator
\begin{eqnarray}\label{Kle}
    \hat\K_{\text{BE}}&=&\sum_{k=1}^K\[\(n-{c_k}-1\)\S_{{k}}^{-1}\]^0-\sum_{k=2}^K\[\(n-{s_k}-1\)\S_{[k]}^{-1}\]^0-\frac{1}{\Tr{\S}}\I
\end{eqnarray}
dominates the MVUE in (\ref{Kub}) in terms of MSE:
\begin{eqnarray}
  E\left\{\left\|\hat\K_{\text{BE}}-\K\right\|^2\right\}\leq E\left\{\left\|\hat\K_{\text{MVUE}}-\K\right\|^2\right\},
\end{eqnarray}
for all  $\K$ in the set ${\mathcal{K}}$ defined in (\ref{KK}).
\end{theorem}
\begin{proof}
 The difference in MSEs is
\begin{eqnarray}\label{deltaml}
  \delta&=&E\left\{\left\|\hat\K_{\text{BE}}-\K\right\|^2\right\}-E\left\{\left\|\hat\K_{\text{MVUE}}-\K\right\|^2\right\}\nonumber\\
  &=&E\left\{-\frac{2\Tr{\hat\K_{\text{MVUE}}}}{\Tr{\S}}+\frac{p}{\Trr{\S}}+2\Tr{\frac{1}{\Tr{\S}}\I\K}\right\}\nonumber\\
  &=&E\left\{-\frac{2\Tr{\hat\K_{\text{MVUE}}}}{\Tr{\S}}+\frac{p}{\Trr{\S}}+2\frac{\Tr{\hat\K_{\text{MVUE}}}}{\Tr{\S}}+4\Tr{\Nabla\frac{1}{\Tr{\S}}\I}\right\}\nonumber\\
  &=&E\left\{\frac{p}{\Trr{\S}}+4\Tr{\Nabla\frac{1}{\Tr{\S}}\I}\right\}\nonumber\\
  &=&E\left\{-\frac{3p}{\Trr{\S}}\right\}\leq 0
\end{eqnarray}
where we applied Theorem \ref{theoremSURE} with
$\H=\frac{1}{\Tr{\S}}\I$ and used the identity
\begin{eqnarray}
  \Tr{\Nabla\frac{1}{\Tr{\S}}\I}=\sum_i \frac{\partial}{\partial \[\S\]_{i,i}} \frac{1}{\Tr{\S}}=-\frac{p}{\Trr{\S}}.
\end{eqnarray}
\end{proof}

Theorem \ref{theorem_scaled} proves the inadmissibility of the MLE and MVUE in any decomposable graphical model. This contribution extends the results in \cite{konno:2001,sunsun:2005,sunsun:2007}. The specific form of $\hat\K_{\text{BE}}$ is not of significant importance and was chosen for simplicity. It is based on a similar Efron-Morris type estimator derived for saturated models in \cite{Tsukuma:2006}.

Finally, it is worth mentioning that Theorem \ref{theorem_scaled} is an example of the well known Stein's phenomenon in which the simultaneous estimation of multiple unrelated parameters can be more accurate than estimating them separately. Indeed, the simplest case of decomposable models is the diagonal (or block diagonal) inverse covariance matrix in which $\S_k$ are statistically independent of each other and depend on different parameters. Theorem \ref{theorem_scaled} establishes that joint estimation can be improved by global shrinkage.

\subsection{SURE based parameter tuning}
The main application of SURE in signal processing is parameter tuning \cite{ulfrasson:08,eldar:gsure,luisier:2007}. Thus, we now illustrate how automatic parameter tuning in decomposable graphical models can utilize SURE.

Consider a class of estimators parameterized by one or more parameters. For simplicity, we restrict ourselves to a special class of estimators with one design parameter:
\begin{eqnarray}\label{Kd}
    \hat\K_d&=&\sum_{k=1}^K\[\(n-{c_k}-1-d\)\S_{{k}}^{-1}\]^0-\sum_{k=2}^K\[\(n-{s_k}-1-d\)\S_{[k]}^{-1}\]^0
\end{eqnarray}
parameterized by $d$. Given this class of estimators, we would like to find the parameter $d$ which minimizes the MSE:
\begin{eqnarray}\label{dex}
  \min_d E\left\{\left\|\hat\K_d-\K\right\|^2\right\},
\end{eqnarray}
or excluding constant terms
\begin{eqnarray}\label{dex1}
  \min_d E\left\{\left\|\hat\K_d\right\|^2-2\Tr{\hat\K_d\K}\right\}.
\end{eqnarray}
Solving (\ref{dex1}) is difficult as the expectation and the second term in the objective depend on $\K$ which is unknown. Instead, we propose to use the SURE result in Theorem \ref{theoremSURE} and replace the unknown MSE with its unbiased estimate:
\begin{eqnarray}
  \min_d \left\|\hat\K_d\right\|^2-2\Tr{\hat\K_d\K_{\text{MVUE}}+2\Nabla\hat\K_d}.
\end{eqnarray}
Substitution of $\hat\K_d$ from (\ref{Kd}) and excluding constant terms yields
\begin{eqnarray}
  \min_d d^2\left\|\D\right\|^2+4d\Tr{\Nabla\D}
\end{eqnarray}
where
\begin{eqnarray}
  \D=\sum_{k=1}^K\[\S_{{k}}^{-1}\]^0-\sum_{k=2}^K\[\S_{[k]}^{-1}\]^0.
\end{eqnarray}
Finally, solving for $d$ results in
\begin{eqnarray}\label{dAB}
  d&=&\frac{-2\Tr{\Nabla\D}}{\left\|\D\right\|^2}\nonumber\\
  &=&\frac{\sum_{k=1}^K \Tr{\S_{{k}}^{-2}}-\sum_{k=2}^K \Tr{\S_{[k]}^{-2}}+\sum_{k=1}^K \Trr{\S_{{k}}^{-1}}-\sum_{k=2}^K \Trr{\S_{[k]}^{-1}}}{\left\|\sum_{k=1}^K\[\S_{{k}}^{-1}\]^0-\sum_{k=2}^K\[\S_{[k]}^{-1}\]^0\right\|^2},
\end{eqnarray}
where we have used (\ref{NablaiS}).

Simulation results presented in Section \ref{sec_num} shows promising performance gains. While we adopted a particularly simple class of estimators in (\ref{Kd}) more advanced estimator classes can likely be treated as well. For example, state-of-the-art methods for covariance estimation in decomposable graphical models involve the use of Bayesian methods and conjugate priors \cite{letac:2007,rajaratman:07}. These distributions depend on tuning parameters that must be chosen beforehand or estimated from the available data. Currently, these parameters are chosen through cross validation, or empirical Bayes methods. It would be interesting to examine the use of SURE as an alternative.

\section{Positive part estimators}\label{sec_pospart}
In the previous sections we proposed estimators which dominate the MLE in terms of MSE. The conditions for their existence are similar to those of the MLE, however they may require more samples in order to be positive semidefinite. For small sample size, we propose to project these estimators onto the set of decomposable positive semi-definite matrices ${\mathcal{K}}$ in (\ref{KK}). We prove that this projection results in legitimate positive semidefinite estimators with better or equal MSE performance.

Let  $\hat\K$ be a given estimator of $\K$. Define $\overline\K$ as the projection of $\hat\K$ onto the set ${\mathcal{K}}$ in (\ref{KK})
\begin{eqnarray}\label{proj}
  \overline\K=\arg\min_{\overline{\K}\in{\mathcal{K}}}\left\|\overline\K-\hat\K\right\|^2.
\end{eqnarray}
The optimization (\ref{proj}) can be expressed as a semidefinite program (SDP). Therefore, the projected estimator $\overline\K$ can be efficiently computed using standard SDP optimization packages, e.g., \cite{sedumi}. The following theorem states that the projected estimator reduces the error with probability one.
\begin{theorem}\label{theorem_proj}
Let  $\hat\K$ be a given estimator of $\K\in {\mathcal{K}}$ and define $\overline\K$ as its projection in (\ref{proj}).
Then,
\begin{eqnarray}\label{wp1}
  \left\|\overline\K-\K\right\|^2\leq \left\|\hat\K-\K\right\|^2,
\end{eqnarray}
with probability one for all $\K$ in the set ${\mathcal{K}}$ in (\ref{KK}).
\end{theorem}
\begin{proof}
The proof is based on the convexity of the set ${\mathcal{K}}$ in (\ref{KK}) and the classical theorem of projection onto convex sets (POCS). POCS states that \cite{boyd:2003}
\begin{eqnarray}
  \Tr{\(\hat\K-\overline\K\)^T\(\K-\overline\K\)}\leq 0,
\end{eqnarray}
for every $\K\in{\mathcal{K}}$. Adding and subtracting $\K$ in the first parenthesis yields
\begin{eqnarray}
  \left\|\K-\overline\K\right\|^2\leq \Tr{\(\K-\hat\K\)^T\(\K-\overline\K\)}.
\end{eqnarray}
Application of the Cauchy Schwartz inequality results in
\begin{eqnarray}
  \left\|\K-\overline\K\right\|^2\leq \left\|\K-\hat\K\right\|\left\|\K-\overline\K\right\|.
\end{eqnarray}
and therefore
\begin{eqnarray}
  \left\|\K-\overline\K\right\|\leq \left\|\K-\hat\K\right\|.
\end{eqnarray}
Since all of the above inequalities apply to any realization of the random matrix $\hat\K$, (\ref{wp1}) holds with probability one.
\end{proof}

When solving (\ref{proj}) is too computationally expensive, we can relax the constraint set and consider the projection onto the semidefinite cone:
\begin{eqnarray}\label{positivepart}
  \min_{\overline{\K}\succeq\0}\left\|\hat\K-\overline\K\right\|^2.
\end{eqnarray}
Similarly to Theorem \ref{theorem_proj}, the semidefinite cone $\left\{\overline{\K}\succeq\0\right\}$ is a convex set and the solution to (\ref{positivepart}) dominates $\hat\K$ in terms of MSE. Its main advantage is that it satisfies a simple closed form. Let
\begin{eqnarray}
  \hat\K=\U\D\U^T
\end{eqnarray}
be the eigenvalue decomposition of $\hat\K$ where $\U$ is a unitary matrix and $\D$ is a diagonal matrix with the eigenvalues $\[\D\]_{ii}$. Then, the projected estimator is equal to
\begin{eqnarray}\label{Kplus}
  \hat\K_+=\U\D_+\U^T
\end{eqnarray}
where $\D_+$ is a diagonal elements with the elements
$\[\D_+\]_{ii}=\max\{\[\D\]_{ii},0\}$. Due to its similarity to
the positive-part shrinkage estimator in James-Stein regression,
we refer to (\ref{Kplus}) as the positive-part estimator.

\section{Numerical results}\label{sec_num}
Here we present results of numerical experiments in order to illustrate the performance of the above estimators. A standard benchmark used for testing (inverse) covariance estimation and covariance selection is the call center data set  \cite{rajaratman:07,bickel:2008,yuan:2007}. Our goal is to demonstrate the estimators precision rather than the model selection accuracy. Therefore, we estimate the true call center covariance matrix using fixed decomposable models as proposed and discussed in \cite{rajaratman:07}. Next, we artificially generate $n$ independent and identically distributed realizations of jointly Gaussian vectors which follow the true covariance structure. We repeat this procedure $100$ times and report the average performance over these independent trials. We use the three decomposable graphical models analyzed in \cite{rajaratman:07}:
\begin{enumerate}
  \item Two coupled cliques: $C_1=\{1,\cdots,70\}$ and $C_2=\{61,\cdots,100\}$.
  \item Banding: an non-stationary autoregressive model with $p=239$ and cliques $C_k=\{j,\cdots,j+L\}$ for $j=1,\cdots,j-p$ with an empirically validated bandwidth of $L=20$.
  \item Differential banding: an empirically validated and refined banding model in which the first $58$ cliques have a bandwidth of $L=14$ and in the following cliques the bandwidth is equal to $L=4$.
\end{enumerate}
Throughout the simulations, we compare the performance of three estimators: the MLE in (\ref{MLEst}), the MVUE in (\ref{Kub}), and the SURE based estimator in (\ref{Kd}) with $d$ given by (\ref{dAB}). At each realization, we compute the estimators and check their semi-definiteness. When an estimator is not positive semidefinite, we resort to its positive part projection defined in (\ref{positivepart}). In Fig. 1-3 we present the normalized MSE defined as $\|\hat{\K}-\K\|^2/\|\K\|^2$ as a function of the sample size $n$.

It is easy to see the significant MSE performance advantage of the MVUE and the SURE based estimators of $\K$ as compared to the MLE. The gain is most significant when the number of samples is small. In this regime, the MLE performs poorly and is actually worse than the zero estimator, i.e., $\hat{\K}=\0$ which ignores the observations altogether, whereas the newly proposed estimators provide reasonable performance. In small sample sizes the MVUE and SURE based estimators had to be adjusted using their positive part variants. Simulation results (not reported) suggest that the improvement in MSE due to the positive part adjustment is negligible.

\begin{figure}\center
\includegraphics[width=0.80\textwidth]{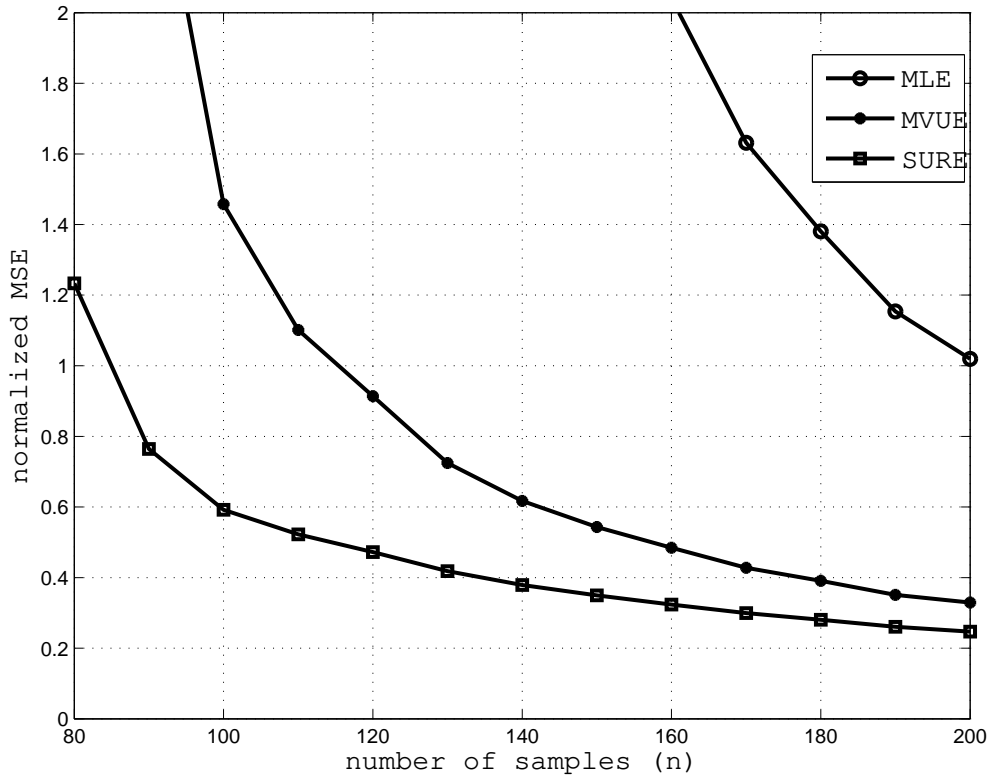}
\caption{Two cliques model: significant MSE improvement with respect to MLE.
\label{fig1}}
\end{figure}

\begin{figure}\center
\includegraphics[width=0.80\textwidth]{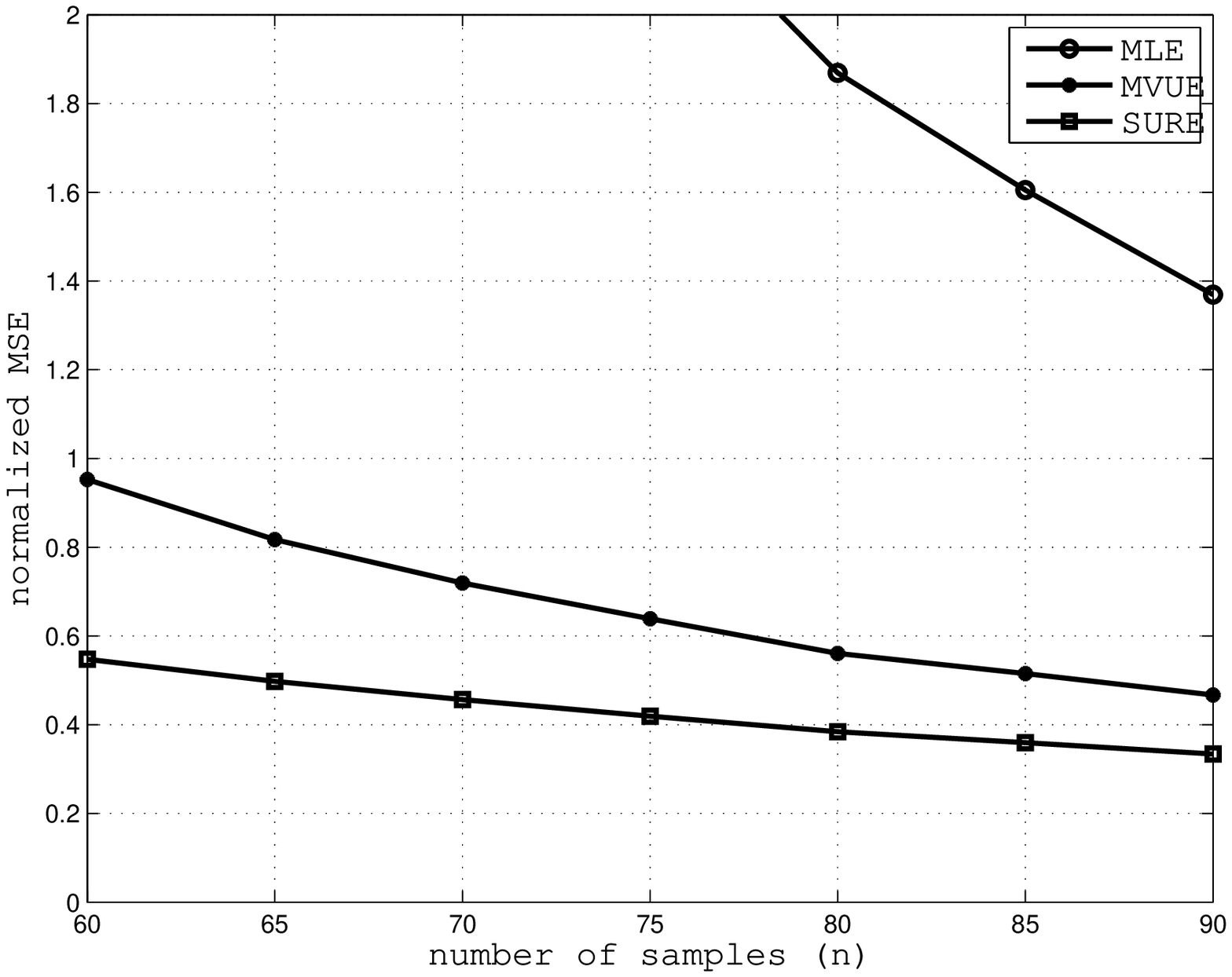}
\caption{Banding model: significant MSE improvement with respect to MLE.
\label{fig2}}
\end{figure}

\begin{figure}\center
\includegraphics[width=0.80\textwidth]{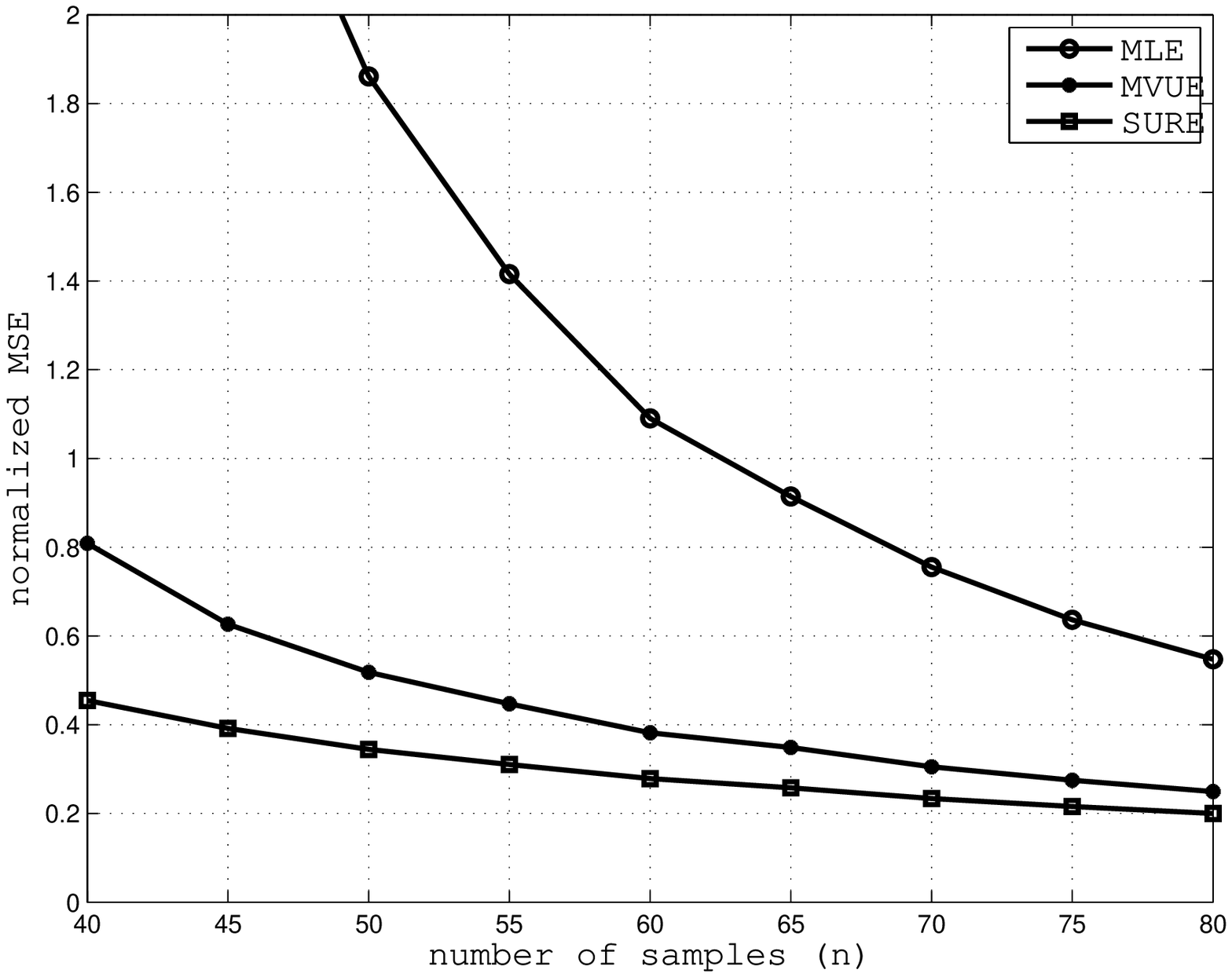}
\caption{Differential banding model: significant MSE improvement with respect to MLE.
\label{fig2}}
\end{figure}

\section{Conclusion and future work}\label{sec_conc}
In this paper, we suggested several alternatives to the MLE for
concentration estimation in decomposable graphical models. We
derived the MVUE and further proposed two biased estimators that
have lower MSE than the MLE. The suggested estimators have simple
closed form solutions and their computational complexity is
similar to the MLE. In addition, we generalized SURE to
decomposable graphical models.

Throughout this work, we assumed that the graphical model is decomposable and illustrated our results for practical signal processing examples, e.g.,  banded and arrow structured concentration matrices. Moreover, any conditional independence graph can be approximated as decomposable using available graph theoretical tools. An important challenge for future work is the extension of our results to non-decomposable graphs.

\section{Acknowledgments}
The authors would like to thank Carlos M. Carvalho for providing the call center dataset.

\appendices
\section{Sparsity in the concentration matrix}\label{app_sparsity}
In this section, we prove that for the Gaussian distribution conditional independence is equivalent to sparsity in the inverse covariance. In order to simplify the notation we define $a=\[\x\]_i$, $b=\[\x\]_j$, and $\c=\[\x\]_{V\setminus i,j}$. The conditional covariance of $a$ and $b$ given $\c$
is defined as
\begin{eqnarray}
  \rho_{ab|\c}=E\left\{\left.\(a-E\left\{\left.a\right| \c\right\}\)
  \(b-E\left\{b\left| \c\right.\right\}\)\right| \c\right\}.
\end{eqnarray}
For the Gaussian distribution this covariance has a well known formula \cite[Theorem 10.2]{kay:93}
\begin{eqnarray}\label{ijrest}
  \rho_{ab|c}=E\left\{ab\right\}-E\left\{a\c^T\right\}
  \[E\left\{\c\c^T\right\}\]^{-1}
  E\left\{\c b\right\}.
\end{eqnarray}
On the other hand, the joint covariance of $a$, $b$, and $\c$ is defined as
\begin{eqnarray}
  \K^{-1}=E\left\{\[\begin{array}{ccc}
                      aa & a\c^T & ab \\ \c a & \c\c^T &\c b \\ ba & b\c^T & bb
                    \end{array}
  \]\right\}.
\end{eqnarray}
Using the matrix inversion lemma for partitioned matrix \cite[page 572]{kay:93}, the top right element of $\K$ is equal to
\begin{eqnarray}
  \[\K\]_{ab}=\Q_{a|c}\rho_{ab|c}\Q_{b|ac},
\end{eqnarray}
where $\rho_{ab|c}$ is the conditional covariance and
\begin{eqnarray}
  \Q_{a|c}&=&\[E\left\{aa\right\}-E\left\{a\c^T\right\}\[E\left\{\c \c^T\right\}\]^{-1}E\left\{\c b\right\}\]^{-1}\\
  \Q_{b|ac}&=&\[E\left\{bb\right\}-E\left\{b\[\begin{array}{cc}
                                                a & \c^T
                                              \end{array}
  \]\right\}\[E\left\{\[\begin{array}{c}
                                                a \\ \c
                                              \end{array}
  \] \[\begin{array}{cc}
                                                a & \c^T
                                              \end{array}
  \]\right\}\]^{-1}E\left\{\c b\right\}\]^{-1}.
\end{eqnarray}
Therefore, if the conditional covariance is equal to zero then the corresponding element in the concentration matrix is also equal to zero:
\begin{eqnarray}
  \rho_{ab|c}=0\quad \Leftrightarrow \quad \[\K\]_{ab}=0.
\end{eqnarray}

\section{MVUE and SURE in the natural exponential family}\label{app exp}
A natural exponential family is defined as
\begin{eqnarray}\label{expdist}
  f\(\x;\thet\)=k\(\x\)e^{\thet^T\x-\psi\(\thet\)}.
\end{eqnarray}
Its natural parameter is $\thet$ and $\x$ is a complete sufficient statistic. The MVUE for estimating $\thet$ given $\x$ is \cite{Hudson:1978}
\begin{eqnarray}\label{expdistMVUE}
  \hat\thet_{\text{MVUE}}=-\frac{\partial }{\partial \x}\log \(k\(\x\)\).
\end{eqnarray}
For any continuously differentiable function of $\x$ denoted by $h_i\(\x\)$, the following SURE identity holds \cite{Hudson:1978,eldar:gsure}
\begin{eqnarray}
  E\left\{h_i\(\x\)\cdot\[\thet\]_i\right\}=E\left\{h_i\(\x\)\[-\frac{\partial \log\(k\(\x\)\)}{\partial \[\x\]_i}\]-\frac{\partial h_i\(\x\)}{\partial
  \[\x\]_i}\right\}.
\end{eqnarray}
Plugging in the MVUE yields
\begin{eqnarray}\label{expdistsure}
  E\left\{h_i\(\x\)\cdot\[\thet\]_i\right\}=E\left\{h_i\(\x\)\[\hat\thet_{\text{MVUE}}\]_i-\Tr{\frac{\partial h_i\(\x\)}{\partial
  \[\x\]_i}}\right\},
\end{eqnarray}
which has an intuitive interpretation: an unbiased estimate of $h_i\(\x\)\[\thet\]_i$ is obtained by simply replacing $\[\thet\]_i$ with its MVUE and adding a correction term.

The Wishart distribution in (\ref{wishart}) belongs to the natural exponential family with parameter $\thet$ being a $\frac{p\(p+1\)}{2}$ vector with the elements in the upper triangular part of $\K$, and variable $\x$ being a  $\frac{p\(p+1\)}{2}$ vector with the elements in the upper triangular part of $\S$ with $-1$ and $-2$ factors on the diagonal and off-diagonal elements, respectively. Therefore, the symmetric matrix version of the MVUE in (\ref{expdistMVUE}) is
\begin{eqnarray}
  \hat\K_{\text{MVUE}}&=&2\overline\Nabla\log \(k\(\S\)\)\nonumber\\
  &=&2\overline\Nabla \log\left|\S\right|^{\frac{n-p-1}{2}}\nonumber\\
  &=&\(n-p-1\)\S^{-1}
\end{eqnarray}
where $\overline\Nabla$ is a symmetric differential operator
\begin{eqnarray}
  \[\overline\Nabla\]_{ij}=\left\{\begin{array}{cc}
                           \frac{\partial }{\partial \[\S\]_{ii}} & i=j \\
                           \frac{1}{2}\frac{\partial }{\partial \[\S\]_{ij}} & i\neq j.
                         \end{array}
  \right.
\end{eqnarray}
This operator compensates for the factor 2 in the off diagonal elements of the derivatives since
\begin{eqnarray}
   \frac{\partial }{\partial \[\S\]_{ij}}\log\left|\S\right|=\[\S^{-1}\]_{ij}+\[\S^{-1}\]_{ji}=2\[\S^{-1}\]_{ij}
\end{eqnarray}
for symmetric matrix $\S$ and $i\neq j$.

Similarly, the SURE expression in (\ref{expdistsure}) is given by \cite{stein:rietz,steinnotes,haff:empirical}
\begin{eqnarray}\label{suremat0}
  E\left\{\Tr{\H\(\S\)\K}\right\}=E\left\{\Tr{\(n-p-1\)\H\(\S\)\S^{-1}+2\overline\Nabla\H\(\S\)}\right\},
\end{eqnarray}
and plugging in the MVUE yields
\begin{eqnarray}\label{suremat}
  E\left\{\Tr{\H\(\S\)\K}\right\}=E\left\{\Tr{\H\(\S\)\hat\K_{\text{MVUE}}+2\overline\Nabla\H\(\S\)}\right\}.
\end{eqnarray}

Gaussian graphical models also belong to the natural exponential family with parameter $\thet$ being a vector with the non-zero elements in the upper triangular part of $\K$, and variable $\x$ being a vector with the corresponding elements of $\S$ and their correction factors. This holds for any Gaussian graphical model, but is not useful unless we can evaluate the function $k\(\cdot\)$ in (\ref{expdist}) and its derivatives. In the case of decomposable models $k\(\cdot\)$ has a simple closed form. Indeed, plugging (\ref{wishart}) into (\ref{hyper}) yields
\begin{eqnarray}
  \log\(k\(\S\)\)=\sum_{k=1}^K\log\left|\S_k\right|^{\frac{n-c_k-1}{2}}-\sum_{k=1}^K\log\left|\S_k\right|^{\frac{n-s_k-1}{2}}.
\end{eqnarray}
Therefore, the MVUE of $\K$ is
 \begin{eqnarray}
    \hat\K_{\text{MVUE}}&=&2\Nabla\log \(k\(\S\)\)\nonumber\\
    &=&\sum_{k=1}^K\[\(n-{c_k}-1\)\S_{{k}}^{-1}\]^0-\sum_{k=2}^K\[\(n-{s_k}-1\)\S_{[k]}^{-1}\]^0,
\end{eqnarray}
and SURE is obtained by modifying the differential operator in (\ref{suremat}) to take into account only the non-zero elements of $\K$ as expressed in (\ref{nabla}).

\section{Technical inequalities}\label{app_proofs}
For simplicity, we partition the submatrix of the $k$th clique as
\begin{eqnarray}
  \S_{{k}}=\[\begin{array}{cc}
          \A & \B \\
          \B^T & \S_{[k]}
        \end{array}
  \].
\end{eqnarray}
Proof of (\ref{CvS}):
Using the partitioned matrix inverse
\begin{eqnarray}
  \[\begin{array}{cc}
          \A & \B \\
          \B^T & \S_{[k]}
        \end{array}
  \]^{-1}=\[\begin{array}{cc}
              \Del^{-1} & -\Del^{-1}\B\S_{[k]}^{-1} \\
              -\S_{[k]}^{-1}\B^T\Del^{-1} & \S_{[k]}^{-1}+\S_{[k]}^{-1}\B^T\Del^{-1}\B\S_{[k]}^{-1}
            \end{array}
  \]
\end{eqnarray}
where $\Del=\A-\B\S_{[k]}^{-1}\B^T$. Therefore,
\begin{eqnarray}
  \Tr{\S_k^{-1}}=\Tr{\Del^{-1}}+\Tr{\S_{[k]}^{-1}}+\Tr{\S_{[k]}^{-1}\B^T\Del^{-1}\B\S_{[k]}^{-1}}\geq \Tr{\S_{[k]}^{-1}},
\end{eqnarray}
where the last inequality is due to the positive semidefiniteness of $\S_k\succeq\0$ and its Schur complement $\Del\succeq\0$.


\noindent Proof of (\ref{CvS2}):
Using the partitioned matrix inverse once again, we obtain
\begin{eqnarray}
  \Tr{\S_k^{-2}}&=&\Tr{\Del^{-2}}+2\Tr{\Del^{-1}\B\S_{[k]}^{-2}\B^T\Del^{-1}}+\Tr{\S_{[k]}^{-2}}\nonumber\\
  &&+2\Tr{\S_{[k]}^{-2}\B^T\Del^{-1}\B\S_{[k]}^{-1}}+\Tr{\(\S_{[k]}^{-1}\B^T\Del^{-1}\B\S_{[k]}^{-1}\)^2}\geq \Tr{\S_{[k]}^{-2}}.
\end{eqnarray}

\end{document}